\documentclass[10pt,twoside]{amsart}
\usepackage{amsmath}
\usepackage{amsthm}
\usepackage{amsfonts}
\usepackage{amssymb}
\usepackage[all]{xy}
\usepackage{alltt}
\usepackage{enumitem}
\usepackage{stmaryrd}
\usepackage{comment}
\usepackage{xspace}
\usepackage{graphicx}
\usepackage{comment}
\usepackage{subcaption}
\theoremstyle{plain}

\theoremstyle{definition}

\usepackage[margin=1.25in]{geometry}

\usepackage[
 colorlinks=true,urlcolor=blue,linkcolor=blue,citecolor=blue
]{hyperref}

\begin{document}  
 
% Article information
\title{A clock model for planetary Conjunctions}
\date{\today}
 
% Author information
\author{Sunil K. Chebolu}
\address{Department of Mathematics \\
llinois State University \\
Normal, IL 61790, USA}
\email{schebol@ilstu.edu}

\maketitle
\thispagestyle{empty}

% Body of the paper

\section{Introduction}\label{sec:introduction}

On December 21, 2020, the night sky offered a beautiful astronomical treat for stargazers worldwide. An exceptionally rare conjunction of Jupiter and Saturn brought them to 0.1 degrees of angular separation -- a fifth of the full moon's diameter. It marked their closest approach since 1623 and the closest visible conjunction since 1226 (almost 800 years ago!). Astronomy enthusiasts crossed their fingers for clear skies and waited eagerly for this event. The internet and social media were inundated with pictures and news reports, celebrating the great conjunction of our solar system's two most extensive and majestic planets. 

The author was fortunate to have clear skies to enjoy this rare celestial event. Below is a picture he took of this conjunction using a Google pixel-3 camera and Celestron 130 mm SLT Newtonian Telescope equipped with a 25 mm eyepiece and a 2x Barlow lens. He saw these two planets in the same telescopic field of view for the first time. 

\begin{figure}[h]
\includegraphics[height=5cm, width=6cm]{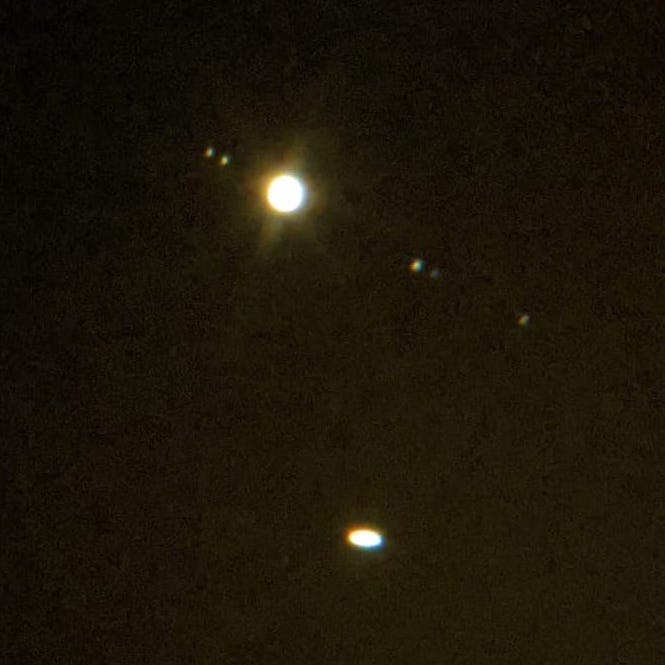}
\caption{Conjunction of Jupiter and Saturn on the winter solstice of 2020}
\end{figure}
 A conjunction of Jupiter and Saturn is called great conjunction because it is the rarest among conjunctions of all pairs of classical planets (Moon, Mercury, Venus, Sun, Mars, Jupiter, and Saturn), happening once every 20 years on average.

This note is written for the mathematicians on the street who are enthusiastic about astronomy. Mathematics students will also appreciate the mathematical framework and analysis for understanding conjunction. We will address the following questions from a mathematical standpoint. What is a conjunction? How often do we have a conjunction of two planets? How are great conjunctions distributed in the sky? How long will it take for a cycle of great conjunctions to return to the same point in the sky and time of the year? Can there be a conjunction of 3 outer planets? How often can that happen? We will begin by setting up the geometric framework for formulating these questions. Then we will analyze these questions by first examining similar questions in a more mundane context: the three hands in a clock. Even though the clock model is a simple model that ignores the more subtle and complex elements of planets' cosmic dance, it captures the main ideas succinctly. We will not discuss the history of great conjunctions. For an excellent historical survey of great conjunctions, see \cite{E}. All numerical calculations were done in the open-source platform SageMath.

\section{Celestial coordinate systems} 

Loosely speaking, two celestial objects (e.g., planets) are in conjunction when they have a close approach from our perspective as one overtakes the other from above or below. We will set the stage for celestial coordinates in the geocentric model to make this precise. These coordinates are similar to the latitudes and longitudes on Earth, allowing us to uniquely identify every location by an ordered pair of numbers.  

\begin{figure}[h]
\includegraphics[height=7cm, width=7cm]{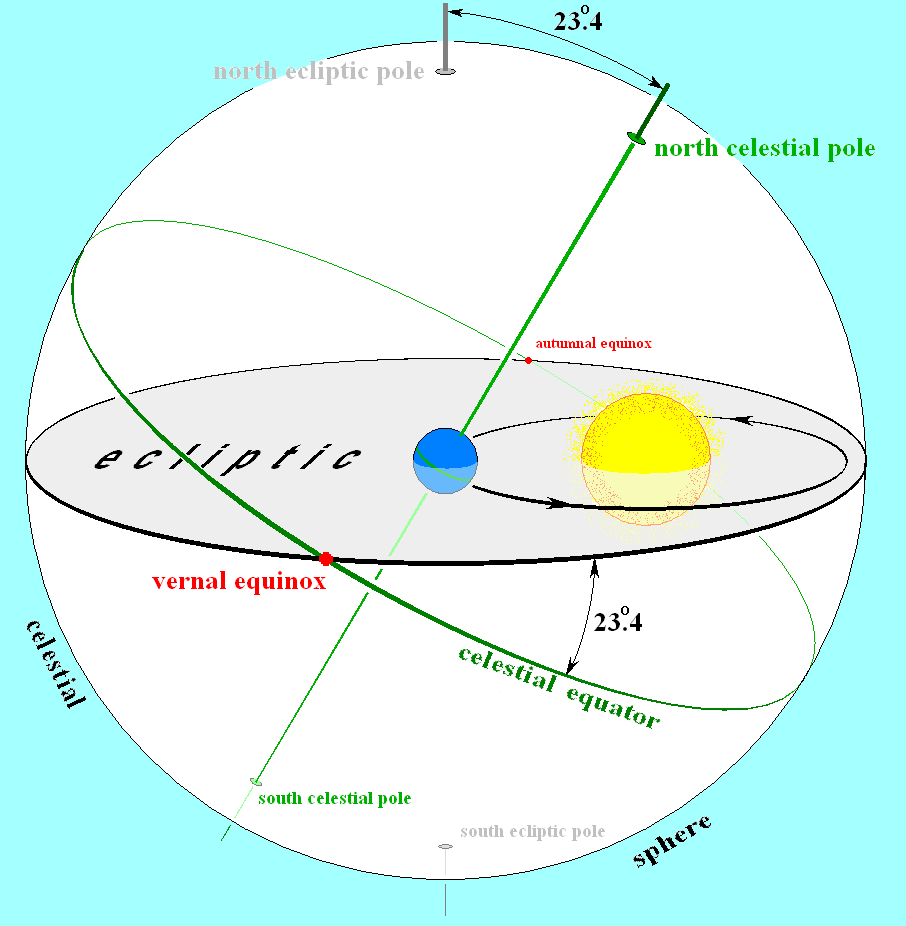}
\caption{The geometry of the celestial sphere. Source: Wikipedia} \label{twoplanes}
\end{figure}

Our Earth will be fixed at the center of the celestial sphere, an imaginary sphere of infinite radius from our vantage point. All celestial objects are projected radially outwards on the celestial sphere. All-stars (except the Sun) are fixed points of light glued onto the celestial sphere that spins on the celestial axis (the axis of Earth's rotation extended both ways to meet the celestial sphere at the celestial poles) once every sidereal day (about 23 hours and 56 minutes). The Sun traces a path on the celestial sphere (relative to the distant stars) that defines the ecliptic plane. Planets orbit in planes that are very close to the ecliptic. The plane that contains the Earth's equator is called the equatorial plane. The Earth's equator projected radially outward on the celestial sphere is the celestial equator. The ecliptic and equatorial planes intersect at an angle of approximately 23.4 degrees. The intersection of the two planes is a line that connects the vernal and autumnal equinoxes going through Earth's center. See figure \ref{twoplanes}.

\begin{figure}[h]
\centering
\begin{subfigure}{.35\textwidth}
 \centering
 \includegraphics[width=.78\linewidth, angle=-28 ]{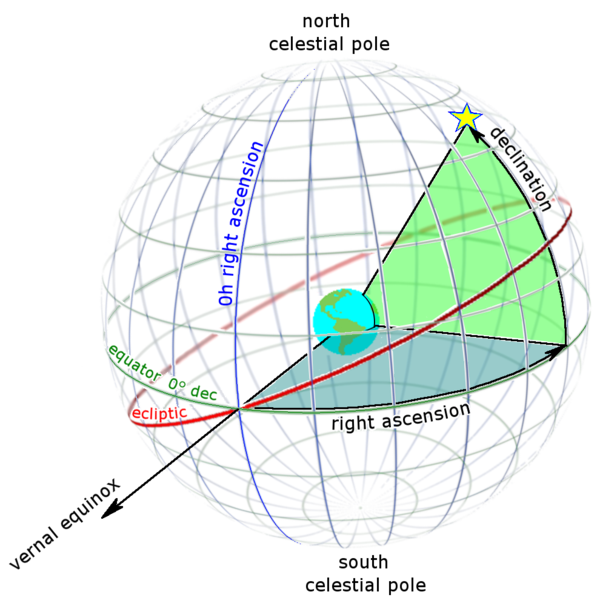}
 \caption{Equatorial Coordinates}
 \label{fig:sub2}
\end{subfigure}
\begin{subfigure}{.5\textwidth}
 \centering
 \includegraphics[width=.7\linewidth]{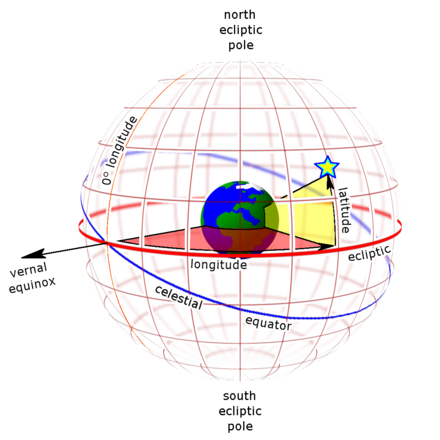}
 \caption{Ecliptic Coordinates}
 \label{fig:sub1}
\end{subfigure}
\caption{Geocentric spherical celestial coordinates. Source: Wikipedia}
\label{2cordinates}
\end{figure}

Recall from geometry that any plane $(P)$ and a vector $\vec{\mathbf{n}}$ that is normal to $P$ define coordinates in the three-dimensional space. Since we will be referring to points on the celestial sphere, it is convenient to talk about spherical coordinates instead of rectangular coordinates. Letting $P$ be the equatorial plane and $\vec{\mathbf{n}}$ be the celestial axis, we get the equatorial coordinate system in which a point on the celestial sphere is determined uniquely by its right ascension and declination, as shown in figure \ref{2cordinates}(A). The object's declination is the angle of elevation from the celestial equator to the object in question. The celestial north pole is at $90^{\circ}$ N; the celestial equator is at $0^{\circ}$. The celestial south pole at $90^{\circ}$ S. The right ascension RA of an object $X$ is measured eastward up to 24h along the celestial equator starting from the first point on Aires (the position of the Sun at the time of vernal equinox); see figure \ref{2cordinates}(A). Similarly, if we let $P$ be the ecliptic plane and $\vec{\mathbf{n}}$ be the ecliptic polar axis, we get the ecliptic coordinates. A point on the celestial sphere is determined by its ecliptic longitude and ecliptic latitude; see figure \ref{2cordinates}(B). 
 For more details see \cite{book}.

We mention a couple of interesting side notes for students of linear algebra. The change of coordinates transformation from ecliptic to equatorial is given by the following rotation matrix ($\phi = 23.4$ degrees)
\[
\begin{bmatrix}
1 & 0 & 0& \\
0 & \cos \phi & \sin \phi \\
0 & -\sin \phi & \cos \phi
\end{bmatrix}
.\]
An eigenvector of this transformation is a non-zero vector pointing to the vernal equinox from Earth's center. 

\section{Conjunction}
\begin{figure}[h]
\centering
 
\begin{subfigure}{.45\textwidth}
 \centering
 \includegraphics[width=1.0\linewidth ]{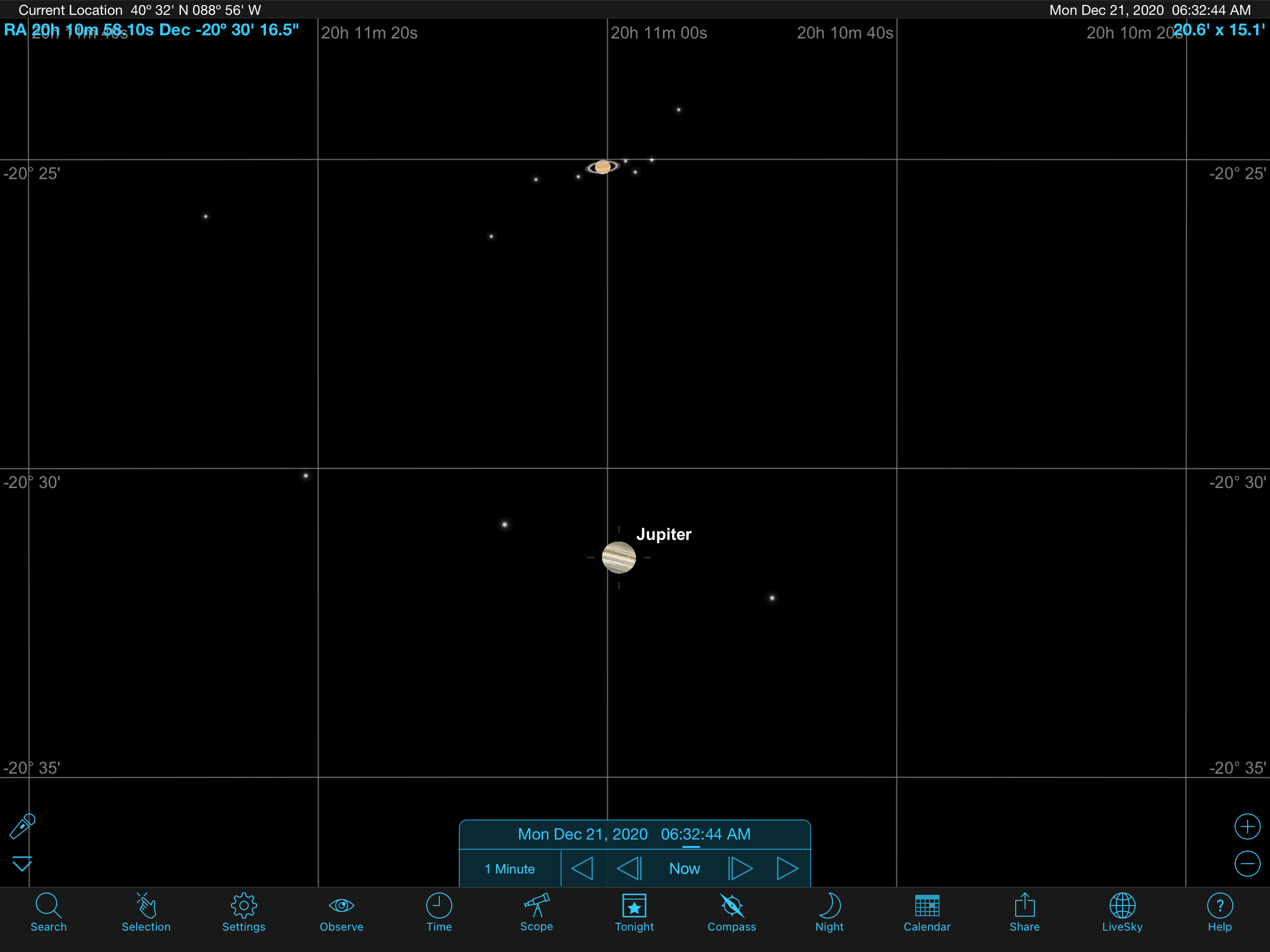}
 \caption{Equatorial Conjunction}
 \label{fig:sub2}
\end{subfigure}
\begin{subfigure}{.45\textwidth}
 \centering
 \includegraphics[width=1.0\linewidth]{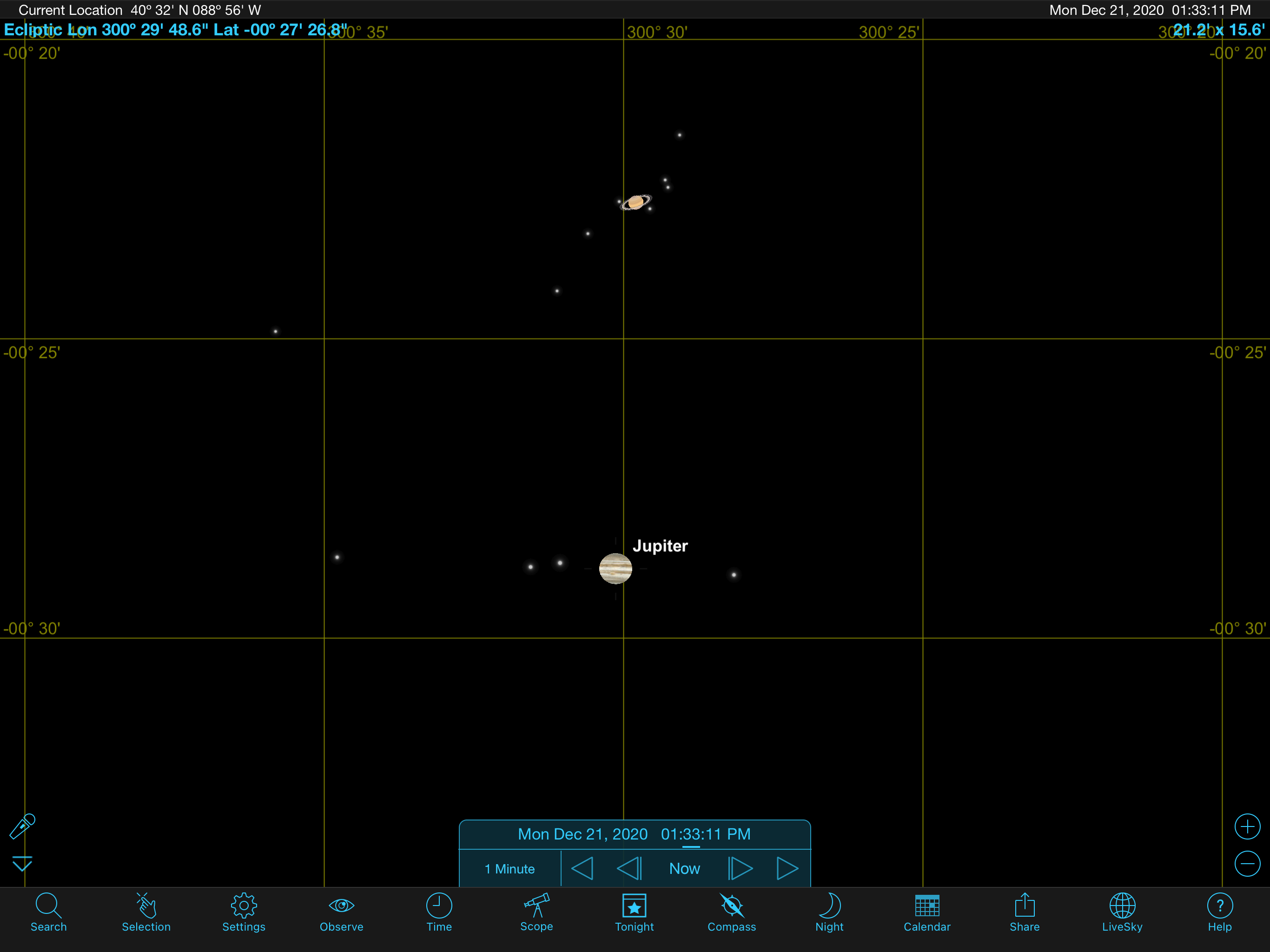}
 \caption{Ecliptic Conjunction}
 \label{fig:sub1}
\end{subfigure}
\caption{Dec 21st, 2020 Equatorial and Ecliptic Conjunctions of Saturn and Jupiter (Source: SkySafari App)}
\label{2020conjunctions}
\end{figure}

Two celestial objects $A$ and $B$ (often planets) are said to be in (geocentric) conjunction if they have either the same right ascension (equatorial conjunction) or the same ecliptic longitude (ecliptic conjunction). \footnote{One can also talk about heliocentric conjunction or conjunction from a third object. In this paper, conjunctions are always in the geocentric sense.}
Note that objects in conjunction need not be close to each other in space or even on the celestial sphere. For instance, in the great conjunction of 2020, Jupiter and Saturn were separated by more than 700 million km. Since planets orbit close to the ecliptic plane, they will appear relatively close during conjunction on the celestial sphere. That is why planetary conjunction is particularly interesting. 

During conjunction, planets may slide from ecliptic conjunction to equatorial or vice versa. For example, Jupiter and Saturn were in ecliptic conjunction on May 28, 2000, and 3 days later, they were in equatorial conjunction. The transition happened the other way around in 2020. On December 21, 2020, Jupiter and Saturn were in equatorial conjunction (RA: 20h 10m 58s) at around 6:30 am CST, and about 7 hours later, they were in ecliptic conjunction (ecliptic longitude: $300^{\circ} 26' 17''$). See figure \ref{2020conjunctions}. Since planets orbit close to the ecliptic plane, a conjunction of planets would most often refer to an ecliptic conjunction.

The next question is, how often does conjunction happen between two planets?

\section{A clock model of conjunctions}
To understand this problem, let us look at a standard 12-hour wall clock - a simple model for the astrometry of planets. Recall that the three hands of this clock spin at different rates. The second-hand takes one minute to complete one circle, the minute-hand takes one hour, and the hour-hand takes 12 hours. When two hands overlap on the top of each other as one passes the other, we say that the two hands are in conjunction. Placing the Sun (or Earth in the geocentric model) at the clock's center, the second, minute, and hour hands correspond respectively (on some scale) to Mars, Jupiter, and Saturn.

\begin{figure}[h]
\includegraphics[height=5cm, width=6.5cm]{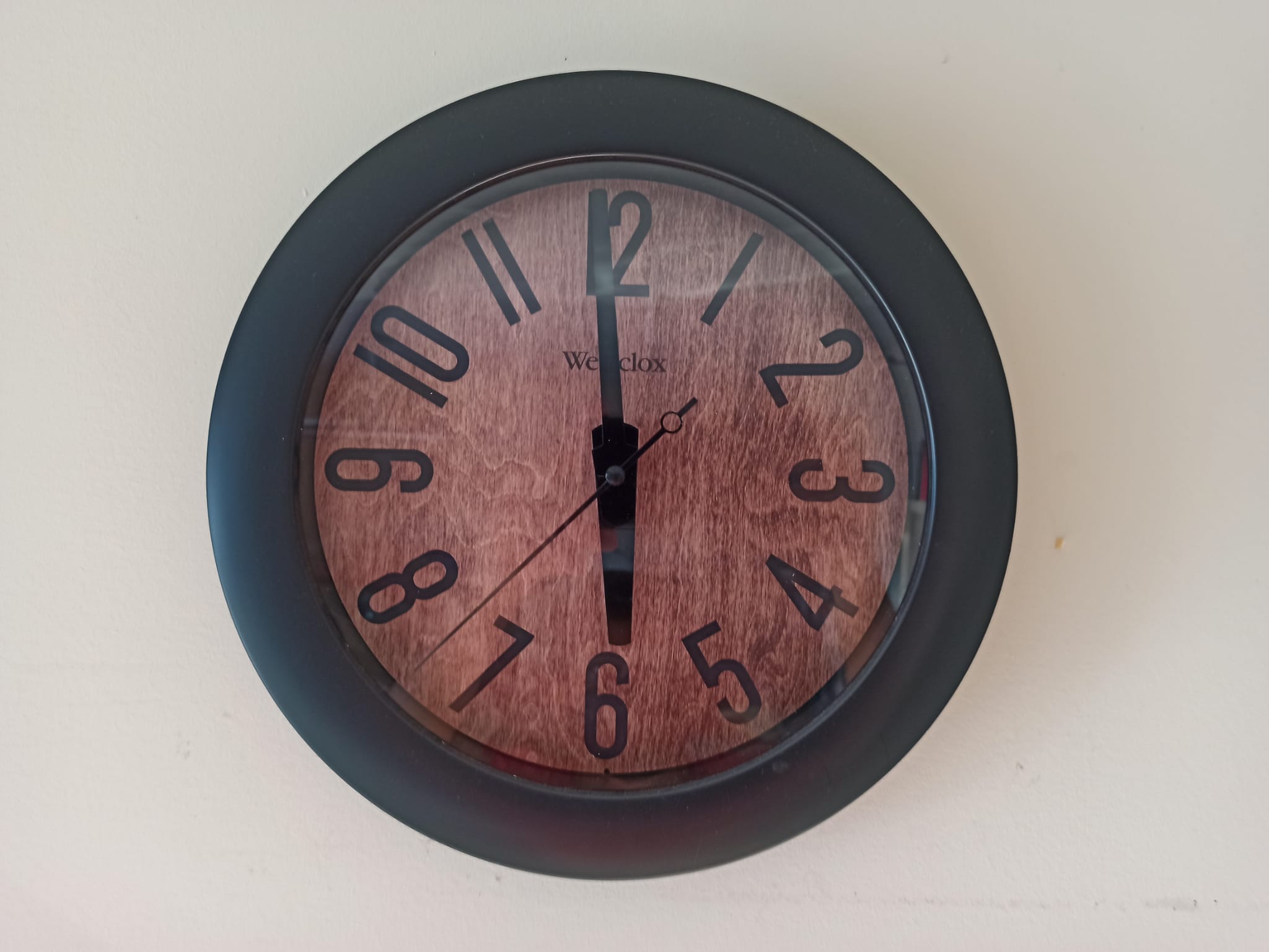}
\caption{The 3 hands of a clock (second, minute, and hour) corresponds to 3 outer planets (Mars, Jupiter, and Saturn).}
\end{figure}

\emph{How often does the minute hand and hours hand come in conjunction?} Let us start with noon when both hands are in conjunction. When does the next conjunction happen? It should happen shortly after 1:05 pm. But exactly when? The hour hand's angular velocity is $\omega_H = 360/12 = 30$ degrees/hr, and the minute hand's angular velocity is $\omega_M = 360/1 = 360$ degrees/hr. The angular distances covered by these hands in $t$ hours is given by
 the functions
$\theta_H(t) = 30t,
\theta_M(t) = 360t
$

Conjunction happens precisely when these two values are equal, modulo $360^{\circ}$. Therefore, finding the moments of conjunctions amounts to solving the equation $\theta_h(t) \equiv \theta_M(t) \mod{360^{\circ}}$:
 \begin{eqnarray*}
30t & \equiv & 360t \mod{360^{\circ}} \\
0 & \equiv & 330t \mod{360^{\circ}}
\end{eqnarray*}
This means conjunction happens at multiplies of $360/330 = 12/11$ hours. So, $t = k\frac{12}{11}$ where $k$ is any non-negative integer. 

In 12 hours, a gap of 12/11 hours ( 1 hr $5.\bar{45}$ minutes) between successive conjunctions gives exactly 11 conjunctions of the hour hand and the minute hand. These conjunctions are distributed on the vertices of a regular hendecagon (11-gon) on the clock with the angle between successive conjunctions equal to $360/11 = 32 \frac{8}{11} ^{\circ}$.

\begin{figure}[h]
\includegraphics[height=5cm, width=5cm]{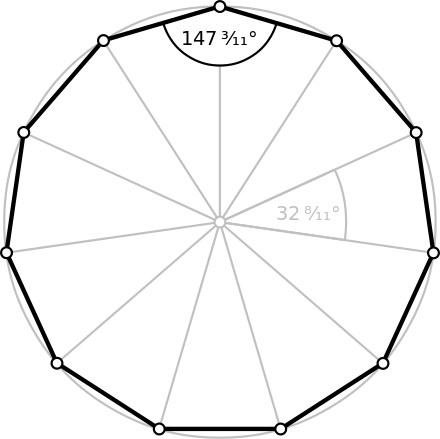}
\caption{The 11 conjunctions of the hour hand and minute hand in 12 hours form the vertices of a regular 11-gon. Source: Wikipedia}
\end{figure}

\section{Jupiter-Saturn (ecliptic) Conjunction}

Planetary conjunction is more complicated than the conjunction of hands of a clock. Unlike the clock's hands, planets do not move in circular orbits (their paths are elliptical). Their speeds are not constant (vary according to Kepler's 2nd laws), and their orbits are tilted to the ecliptic plane. Nevertheless, ignoring these technicalities, the clock model still captures the main idea and helps derive the average distribution of conjunctions in time and space.

We begin by recording a general result. Consider objects $A$ and $B$ that move along a circle at a uniform speed with time periods $t_A$ and $t_B$ (assume $t_A \le t_B$). Then their angular velocities are $\omega_A = 360/t_A$ and $\omega_B = 360/t_B$. Proceeding as in the analysis of hours and minutes hand, we see that conjunction happens at multiples of $360/(\omega_A - \omega_B) = t_A t_B/( t_B - t_A)$.  

Planets orbit around the Sun, but our celestial coordinates are earth-centered, so we must look at planets' angular velocities relative to Earth. The angular velocity of a planet ($P$) relative to Earth ($E$) is given by $\omega_{PE} := \omega_P - \omega_E$. Jupiter takes approximately 12 earth-years, and Saturn takes approximately 30 earth-years to orbit the Sun. Since $12 \times 5 = 30 \times 2$, we have a $5:2$ resonance between the two planets. That is, 5-Jupiter years is equal to 2-Saturn years. The angular velocities are given by
\[ \omega_E = 360/1 = 360, \omega_J = 360/12 = 30, \omega_S = 360/30 = 12,\]
degrees per year.
The average time between successive (ecliptic) conjunctions of Jupiter and Saturn is 
\[360/(\omega_{SE} - \omega_{JE}) = 360/(\omega_{S} - \omega_{J}) = 360/(30 -12) = 360/18 = 20 \, \text{years}.\]
So, on average, great conjunction happens every 20 years. 

In 20 years, Saturn completes about $2/3$ rd of its orbit around the Sun. So successive great conjunctions are separated by $240^{\circ}$. Therefore the three successive great conjunctions are approximately $120^{\circ}$ apart, forming the vertices of an equilateral triangle on the ecliptic equator.

 However, taking more accurate orbital periods	(11.86 years for Jupiter and 29.46 years for Saturn) reveals a different picture. The average period between great conjunctions is closer to 19.85 years (using the above formula). The radial vector from the Sun to Saturn swings approximately $(360/29.46)19.85 = 245.56^{\circ} $ between two consecutive great conjunctions. 
 
 Recall that any point on the plane can be expressed using cartesian coordinates $(x, y)$ or a complex number ($x+iy$). For our purposes, complex numbers are well-suited because rotation by angle $\theta$ of a point about the origin is given by multiplication by $e^{i \theta}$.  

\begin{figure}[h]
\includegraphics[height=11cm, width=12cm]{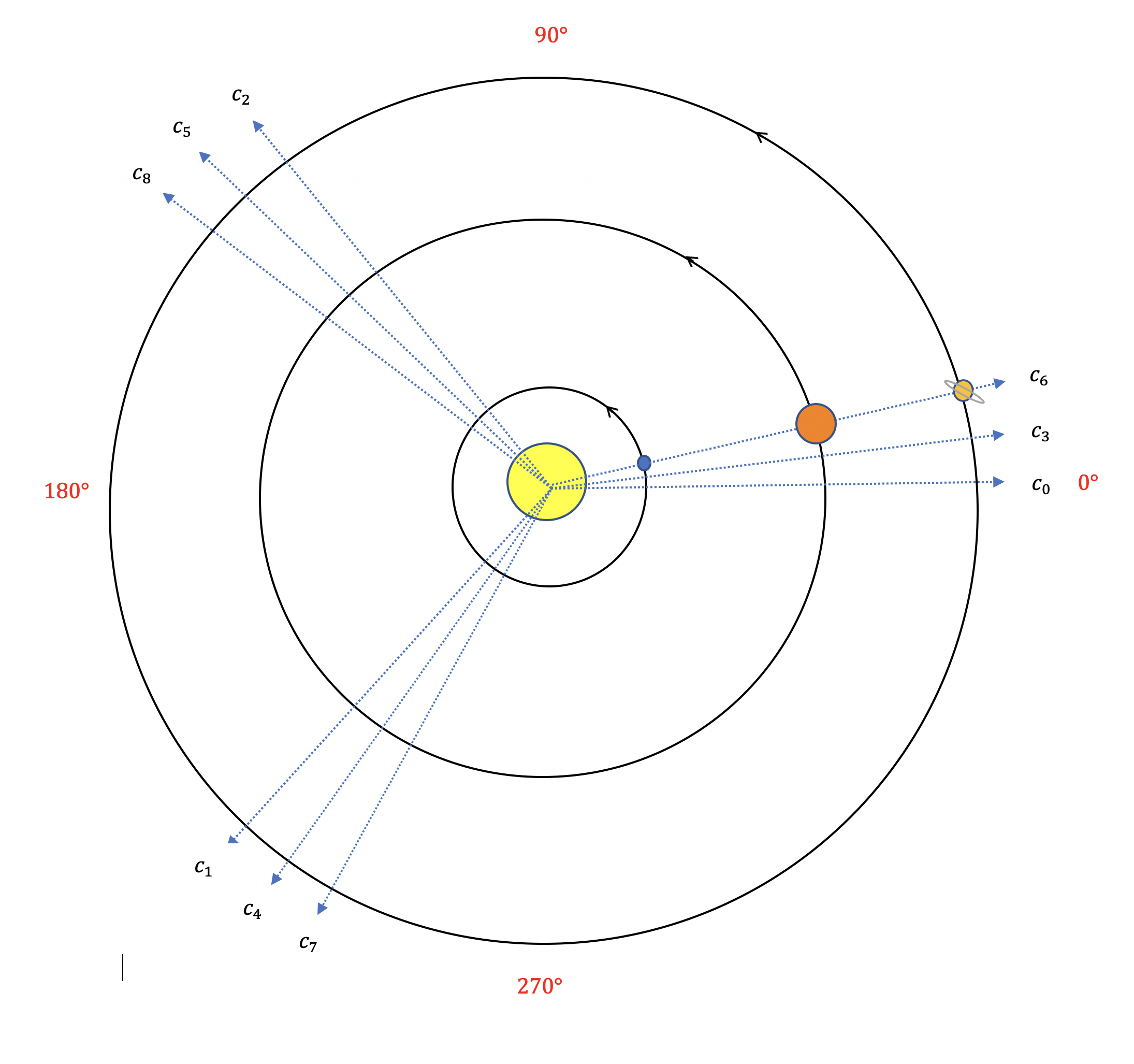}
\caption{9 consecutive Saturn-Jupiter conjunctions (figure not drawn to scale)}
\label{Kepler}
\end{figure}

In this framework of complex numbers, the successive great conjunctions are given by the complex numbers
\[ c_n := re^{i245.56 n} = r (\cos(245.56 n) + i \sin(245.56 n) )\]
on the complex plane where the Sun is at the center, and $r$ is the mean distance between Sun and Saturn. Figure \ref{Kepler} shows a plot of these complex numbers. Let $C_0$ be one of the great conjunctions (like in December 2020). The great conjunctions $C_n$ form 3 different families: $n \equiv 0 \mod 3$, $n \equiv 1 \mod 3$, $n \equiv 2\mod 3$. In particular, note that every 4th great conjunction does not return to the 1st great conjunction; it happens 16.68 degrees east of the 1st conjunction. 

The triangle $\Delta_n$ formed by the vertices $C_{3n}, C_{3n+1},$ and $C_{3n+2}$ is rotating in the same direction as the planets in small ticks (for each tick, the triangle spins 16.68 degrees in 59.55 years). Kepler first observed this fact. 

\begin{center}
\begin{tabular}{ |c|c|c|c| } 
\hline
$C_n$ & degrees from $C_0$&number of Earth years \\
\hline
$C_0$ & $0^{\circ}$ & 0 \\
$C_1$ & $245.56^{\circ}$ &19.85\\
$C_2$ &$131.12^{\circ}$ &39.70\\
$C_3$ &$16.68^{\circ}$ & 59.55\\
$C_4$&$262.24^{\circ}$ &79.40\\
$C_5$ &$147.80^{\circ}$ &99.25\\
$C_6$& $33.36^{\circ}$&119.10\\
$C_7$& $278.92^{\circ}$ &138.95\\
$C_8$&$164.48^{\circ}$&158.80 \\
\hline
\end{tabular}
\end{center}

The table shows nine consecutive great conjunctions. It gives their positions ($245.56n \mod{360^{\circ}}$, measured counter-clockwise from the $C_0$ conjunction) and the number of Earth-years lapsed ($19.85n$) since $C_0$. To find a cycle of great conjunctions that repeats in space (based on celestial longitudes) and in time (based on the solar calendar), we need to find an integer $k$ with the following two properties.
\begin{enumerate}
\item $245.56k \mod{360^{\circ}}$ is as close as possible to $0^{\circ}$: This will bring the cycle of the great conjunction of length $k$ in sync with the ecliptic coordinates. As observed above, $k=3$ offsets by $16.67^{\circ}$.\\
\item $19.85k$ is as close as possible to an integer: This will bring the cycle of the conjunction of length $k$ in sync with the solar calendar. In the year's column, notice that the first number closest to a whole number is 119.10. The Earth returns to approximately the same place in its orbit in 119.10 years. So every sixth great conjunction will look near identical to the first one.\\
\end{enumerate}

The length of the cycles increases if one demands greater precision in the synchronization. For instance, taking more rows of the above table, we see that $C_{22}$ comes closest to $C_0$ great conjunction to within $2.31$ degrees in $436.7$ years- 4 centuries! $C_{66}$ comes close to $C_0$ within 6.96 years in 1310 years - 13 centuries! 

Not all conjunctions are picturesque. First, some of these great conjunctions happen too close to the Sun, making it difficult or impossible to see at most latitudes. Second, the angular separation during great conjunction can vary between $0.1^{\circ}$ to $1.3^{\circ}$ (that is one-fifth to 2.5 times the width of the full moon). The main reason for this variation is the different inclinations of the orbital planes of these planets with the ecliptic (Jupiter inclines $1.31^{\circ}$ and Saturn $2.49^{\circ}$) and the position of the conjunction. The orbital paths meet the ecliptic at two points called the nodes. These nodes are roughly the same for both Saturn and Jupiter. 
The closest great conjunctions happen when the two planets are closer to these nodes because the planets' angular separation will be minimal at the nodes. That is what made the December 2020 great conjunction so special. The two planets came very close (0.1 degrees apart), and it happened $30^\circ$ east of the Sun. So it was a beautiful spectacle for those who had clear skies.

\section{Mars-Jupiter-Saturn Conjunction}

Now let us add Mars to this cosmic dance. How often can Mars, Saturn, and Jupiter be in conjunction? 

We can return to the clock model and ask a similar question. How often will there be a conjunction of the second hand, minutes hand, and hour hand? We already saw that the conjunction between the hour hand and minute hand happens every $12/11$ hours.
Similar, we can check using our formula that a conjunction between a second hand and minute hand occurs every $1/59$ hours:
\[ t_St_M/(t_M-t_S) = 1/60/(1-1/60) = 1/59 \ \ \text{hours}\]

All three hands of the clock will be in conjunction every $t$ hours, where $t$ is the least common integer multiple (lcm) of $12/11$ and $1/59$ hours. Let us recall a number-theoretic fact.

Given two positive rational numbers $a/b$ and $c/d$, it can be shown that 
\[\text{lcm} \left(\frac{a}{b}, \frac{c}{d} \right) = \frac{\text{lcm}(a, c)}{\text{gcd}(b, d)}.\]
Therefore $\text{lcm}(12/11, 1/59) = \text{lcm}(12, 1) = 12$, because $11$ and $59$ are primes.

That means the clock's three hands will be in conjunction once every 12 hours, or twice a day - at noon and midnight.

Perfect conjunction between two planets orbiting the Sun at different angular speeds always happens when the faster-moving planet overtakes the slower planet. However, perfect conjunction between more than two planets is a rarer event. Strictly speaking, a perfect conjunction of more than two objects seldom happens. For Mars, Jupiter, and Saturn to be in conjunction, the pairs (Mars, Jupiter) and (Jupiter, Saturn) have to be in conjunction simultaneously. 

Let $t_M, t_J$, and $t_S$ denote the orbital periods of Mars, Jupiter, and Saturn, respectively.
 We may assume that these are rational numbers. Mars and Jupiter will be in conjunction once every ${t_Jt_M}/{(t_J-t_M)}$ years, and Jupiter and Saturn once every ${t_St_J}/{(t_S-t_J)}$. 
 So the the average time period between successive triple conjunctions is given by 
\[ \text{lcm} \left( \frac{t_St_J}{(t_S-t_J)}, \frac{t_Jt_M}{(t_J-t_M)} \right).\]

Setting $t_M= 1.8, t_J= 12$ and $t_S= 30$, we have the following. Conjunction between Mars and Jupiter happens every 36/17 years, between Jupiter and Saturn every 20 years, and between Mars and Saturn every 90/47. The lcm of any two of these pairs is 180 years. So, on average, once every 180 years, we have this triple conjunction.

The calculation is sensitive to the degree of approximation of the orbital periods. For instance, setting $t_S= 29.5$ and $t_J= 11.8$ years (more accurate values), it turns out that triple conjunction occurs every 531 years, on average!

Great conjunctions are beautiful treats. The next close visible great conjunction will be on March 15, 2080, when 6.0 arcminutes separate Jupiter and Saturn. Mark your calendars and ask your children and grandchildren to enjoy this event.

%https://sagecell.sagemath.org/?q=ajtezd

%https://sparky.rice.edu/public-night/jupsat2.html

%\begin{figure}[h]
%\includegraphics[height=6cm, width=4cm]{palm.jpg}
%\caption{12 knuckles $\implies$12 hours of day/night}
%\end{figure}

\bibliographystyle{plain}

\end{document}